\documentclass{amsart}
\usepackage{amsmath,amssymb,amsfonts}
\usepackage{eucal}
\usepackage{verbatim}

\usepackage{enumerate}

\numberwithin{equation}{section}

\theoremstyle{plain}
\newtheorem{theorem}[equation]{Theorem}

\newtheorem{prop}[equation]{Proposition}

\theoremstyle{definition}
\newtheorem{definition}[equation]{Definition}

\newtheorem*{thank}{Acknowledgments}

\input xy
\xyoption{all}

\newcommand{\Deltaop}{{\bf \Delta}^{op}}

\newcommand{\Hom}{\text{Hom}}
\newcommand{\SSets}{\mathcal{SS}ets}

\newcommand{\Sets}{\mathcal Sets}

\newcommand{\Secat}{\mathcal Se \mathcal Cat}
\newcommand{\Ho}{\text{Ho}}
\newcommand{\map}{\text{map}}
\newcommand{\Map}{\text{Map}}
\newcommand{\ob}{\text{ob}}
\newcommand{\hoequiv}{\text{hoequiv}}

\newcommand{\css}{\mathcal{CSS}}

\newcommand{\Qcat}{\mathcal{QC}at}
\newcommand{\Cat}{\mathcal Cat}

\begin{document}

\title{A Survey of $(\infty,1)$-Categories}

\author[J.E. Bergner]{Julia E. Bergner}

\address{Department of Mathematics \\
Kansas State University \\
Manhattan, KS 66506}

\email{bergnerj@member.ams.org}

\date{\today}


\begin{abstract}
In this paper we give a summary of the comparisons between
different definitions of so-called $(\infty,1)$-categories, which
are considered to be models for $\infty$-categories whose
$n$-morphisms are all invertible for $n>1$.  They are also, from
the viewpoint of homotopy theory, models for the homotopy theory
of homotopy theories.  The four different structures, all of which
are equivalent, are simplicial categories, Segal categories,
complete Segal spaces, and quasi-categories.
\end{abstract}

\maketitle

\section{Introduction}
The intent of this paper is to summarize some of the progress that
has been made since the IMA workshop on $n$-categories on the
topic of $(\infty,1)$-categories.  Heuristically, a
$(\infty,1)$-category is a weak $\infty$-category in which the
$n$-morphisms are all invertible for $n>1$. Practically speaking,
there are several ways in which one could encode this information.
In fact, at the moment, there are four models for
$(\infty,1)$-categories: simplicial categories, Segal categories,
complete Segal spaces, and quasi-categories. They have arisen out
of different motivations in both category theory and homotopy
theory, but work by the author and Joyal-Tierney has shown that
they are all equivalent to one another, in that they can be
connected by chains of Quillen equivalences of model categories.

From the viewpoint of higher category theory, these comparisons
provide a kind of baby version of the comparisons which are being
attempted between various definitions of weak $n$-category.  In
\cite{toen2}, To\"{e}n actually axiomatizes a theory of
$(\infty,1)$-categories and proves that any such theory is
equivalent to the theory of complete Segal spaces.  In
\cite{toen1}, he sketches arguments for proving the equivalences
between the four structures used in this paper.  Although some of
the functors he suggests do not appear to give the desired Quillen
equivalences (or at any rate are not used in the known proofs), he
gives a good overview of the problem.  Another good introduction
of the problem can be found in a preprint by Porter \cite{porter},
and a nice description of the idea behind $(\infty,1)$-categories
can be found in \cite{toenstack}.

From another point of view, these comparisons are of interest in
homotopy theory, as what we are here calling a
$(\infty,1)$-category can be considered to be a model for the
homotopy theory of homotopy theories, a concept which will be made
more precise in the section on simplicial categories.  The idea is
that a simplicial category is in some way ``naturally" a homotopy
theory but functors between such are not particularly easy to work
with.  The goal of finding an equivalent but nicer model led to
Rezk's complete Segal spaces \cite{rezk} and the task of showing
that they were essentially the same as simplicial categories.

It should be noted that there are other proposed models, including
that of $A_\infty$-categories.  Joyal and Tierney briefly discuss
some other approaches to this idea in their epilogue \cite{jt}.

Furthermore, we should also mention that these structures are of
interest in areas beyond homotopy theory and higher category
theory.  For instance, there are situations in algebraic geometry
in which simplicial categories have been shown to provide
information that the more commonly used derived category cannot.
Given a scheme $X$, for example, its derived category $\mathcal
D(X)$ does not seem to determine its $K$-theory spectrum, whereas
its simplicial localization $\mathcal L(X)$ does \cite{tv}.
Furthermore, the simplicial category $\mathcal L(X)$ forms a
stack, which $\mathcal D(X)$ does not \cite{hs}.  Similar work is
also being done using dg categories, which are in many ways
analogous to simplicial categories \cite{toendg}.  In particular,
the category of dg categories has a model category structure which
is defined using the same essential ideas as the model structure
on the category of simplicial categories \cite{tab}.

Also motivated by ideas in algebraic geometry, Lurie uses
quasi-categories and their relationship with simplicial categories
in his work on higher stacks \cite{lurie}.  The first chapter of
his manuscript is also a good introduction to many of the ideas of
$(\infty,1)$-categories.

Another application of the model category of simplicial categories
can be found in recent work of Douglas on twisted parametrized
stable homotopy theory \cite{doug}.  He uses diagrams of
simplicial categories weakly equivalent to the simplicial
localization of the category of spectra (i.e., equivalent as
homotopy theories to the homotopy theory of spectra) in order to
define an appropriate setting in which to do Floer homotopy
theory.

In this paper, we will provide some background on each of the four
structures and then describe the various Quillen equivalences.

\begin{thank}
I would like to thank Bill Dwyer, Chris Douglas, Andr\'e Joyal,
Jacob Lurie, Peter May, and Bertrand To\"en for reading early
drafts of this paper, making suggestions, and sharing their work
in this area.
\end{thank}

\section{Background on model categories and simplicial sets}

Since this paper is meant to be an overview, we are not going to
go deeply into the details here, and the reader familiar with
model categories and simplicial sets may skip to the next section.
For the non-expert, we will give the basic ideas behind both model
category structures and simplicial sets, as well as other
simplicial objects.

The motivation for a model category structure is a common
occurrence in many areas of mathematics.  Suppose that we have a
category whose morphisms include some, called weak equivalences,
which we would like to think of as isomorphisms but do not
necessarily have inverses. Two classical examples are the category
of topological spaces and (weak) homotopy equivalences between
them, and the category of chain complexes of modules over a ring
$R$ and the quasi-isomorphisms, or morphisms which induce
isomorphisms on all homology groups. In order to make these maps
actually isomorphisms, we could formally invert them. Namely, we
could formally add in an inverse to every such map and then add
all the necessary composites so that the result would actually be
a category.  The problem with this approach is that often this
process results in a category such that the morphisms between any
two given objects form a proper class rather than a set. While it
is common to work in categories with a proper class of objects, it
is generally assumed that there is only a set of morphisms between
any two objects, even if there is a proper class of morphisms
altogether.

Imposing the structure of a model category allows us formally to
invert the weak equivalences while keeping the morphisms under
control.  A model category $\mathcal M$ is a category with three
distinguished classes of morphisms, the weak equivalences as
already described, plus fibrations and cofibrations, satisfying
several axioms.  We refer the reader to \cite{ds}, \cite{hirsch},
\cite{hovey}, or the original \cite{quillen} for these axioms. The
importance of these axioms is that they allow us to work with
particularly nice objects in the category, called
fibrant-cofibrant objects, between which we can define ``homotopy
classes of maps" even in a situation where the traditional notion
of homotopy class (as in topological spaces) no longer makes
sense.  The axioms of a model category guarantee that every object
of $\mathcal M$ has a fibrant-cofibrant replacement, and thus we
can define the \emph{homotopy category} of $\mathcal M$, denoted
$\Ho(\mathcal M)$, to be the category whose objects are the same
as those of $\mathcal M$ and whose morphisms are homotopy classes
of maps between the respective fibrant-cofibrant replacements.  In
fact, this construction is independent of the choice of such
replacements, and the homotopy category, up to equivalence, is
independent of the choice of fibrations and cofibrations.  There
are many examples of categories with two (or more) different model
category structures, each leading to the same homotopy category
because the weak equivalences are the same even if the fibrations
and cofibrations are defined differently.

As an example of a model category, consider the category of
topological spaces and the subcategory of weak homotopy
equivalences, or maps which induce isomorphisms on all homotopy
groups.  There is a natural choice of fibrations and cofibrations
such that this category has a model category structure.  (There is
also a model category structure on this category where the weak
equivalences are the homotopy equivalences, but the former is
considered to be the standard model structure.)

One can then define what it means to have a map between model
categories, namely, a functor which preserves essential properties
of the model structures.  It is convenient to use adjoint pairs of
functors to work with model structures, where the left adjoint
preserves cofibrations and the right adjoint preserves fibrations.
Such an adjoint pair is called a Quillen pair. There is also the
notion of Quillen equivalence between model categories, where the
adjoint functors preserve the essential homotopical information.
In particular, a Quillen equivalence induces an equivalence of
homotopy categories, but it is in fact much stronger, in that it
preserves higher-order information, as we will discuss further in
the next section.

Given this kind of structure and interaction between structures,
one can ask the following question, one that can, in fact, be
considered a motivation for the research described in this paper.
Given a model category $\mathcal M$, is there a model category
$\mathcal N$ which is Quillen equivalent to $\mathcal M$ but which
more easily provides information that is difficult to obtain from
$\mathcal M$ itself?  The properties one might look for in
$\mathcal N$ depend very much on the question being asked.  Thus,
it is hoped that the four model structures given in this paper
will be able to provide information about one another.

An important illustration is that of topological spaces and
simplicial sets.  Heuristically, simplicial sets provide a
combinatorial model for topological spaces, and the fact that they
are a ``model" here means that there is a model category structure
on the category of simplicial sets which is Quillen equivalent to
the standard model category structure on the category of
topological spaces.  Simplicial sets are frequently (but not
always) easier to work with because they are just combinatorial
objects.

To give a formal definition of a simplicial set, consider the
category ${\bf \Delta}$ of finite ordered sets $[n]=\{0
\rightarrow 1 \rightarrow 2 \rightarrow \cdots \rightarrow n\}$
for each $n \geq 0$, and order-preserving maps between them.  (As
the notation suggests, one can also consider $n$ to be a small
category.) Let $\Deltaop$ denote the opposite category, where we
reverse the direction of all the morphisms. Then a simplicial set
is a functor $X:\Deltaop \rightarrow \Sets$. There is a geometric
realization functor between the category $\SSets$ of simplicial
sets and the category of topological spaces. Specifically, an
element of $X_0$ is assigned to a point, an element of $X_1$ is
assigned to a geometric 1-simplex, and so forth, where
identifications are given by the face maps of the simplicial set.
Thus, a simplicial set can be regarded as a generalization of a
simplicial complex, where the simplices are not required to form
``triangles" and a given simplex of degree $n$ is regarded as a
degenerate $k$-simplex for each $k>n$ \cite[I.2]{gj}.

In fact, we can perform this kind of construction in categories
other than sets.  A simplicial object in a category $\mathcal C$
is just a functor $X: \Deltaop \rightarrow \mathcal C$.  The
primary example we will consider in this paper is that of
simplicial spaces, or functors $\Deltaop \rightarrow \SSets$.  To
emphasize the fact that they are simplicial objects in the
category of simplicial sets, they are often also called
bisimplicial sets.

Given these main ideas, we can now turn to the four different
models for $(\infty,1)$-categories, or homotopy theories.

\section{Simplicial categories}

The first of the four categories we consider is that of small
simplicial categories.  By a simplicial category, we mean what is
often called a simplicially enriched category, or a category with
a simplicial set of morphisms between any two objects.  Given two
objects $a$ and $b$ in a simplicial category $\mathcal C$, this
simplicial set is denoted $\Map_\mathcal C(a,b)$.  This
terminology is potentially confusing because the term ``simplicial
category" can also be used to describe a simplicial object in the
category of all small categories.  We recover our sense of the
term if the face and degeneracy maps are all required to be the
identity map on objects.

Simplicial categories have been studied for a variety of reasons,
but here we will focus on their importance in homotopy theory,
and, in particular, on how a simplicial category can be considered
to be a homotopy theory.

We should note here that although for set-theoretic reasons we
restrict ourselves to small simplicial categories, or those with
only a set of objects, in practice many of the simplicial
categories one cares about are large.  The standard approach to
this problem is to assume that one is working in a larger universe
for set theory in which the given category is indeed ``small."

Given a model category $\mathcal M$, we can consider its homotopy
category $\Ho(\mathcal M)$.  For many applications, it is
sufficient to work in the homotopy category, but it is important
to remember that in passing from the original model category to
the homotopy category we have lost a good deal of information.  Of
course, part of the goal was formally to invert the weak
equivalences, but in addition the model category possessed
higher-homotopical information that the homotopy category has
lost.  For example, the model category contains the tools needed
to take homotopy limits and homotopy colimits.

In a series of papers, Dwyer and Kan develop the theory of
simplicial localizations, in which, given a model category
$\mathcal M$, one can obtain a simplicial category which still
holds this higher homotopical information.  In fact, they
construct two different such simplicial categories from $\mathcal
M$, the standard simplicial localization $LM$ \cite{dksl} and the
hammock localization $L^H\mathcal M$ \cite{dkfc}, but the two are
equivalent to one another \cite[2.2]{dkcsl}. Furthermore, taking
the component category $\pi_0 LM$, which has the objects of $LM$
and the morphisms the components of the simplicial hom-sets of
$LM$, is equivalent to the homotopy category $\Ho(\mathcal M)$
\cite{dksl}.

Furthermore, there is a natural notion of ``equivalence" of
simplicial categories, which is often called a \emph{Dwyer-Kan
equivalence} or simply \emph{DK-equivalence}.  It is a
generalization of the definition of equivalence of categories to
the simplicial setting.  In particular, a DK-equivalence is a
simplicial functor $f: \mathcal C \rightarrow \mathcal D$ between
two simplicial categories satisfying the following two conditions:
\begin{enumerate}
\item For any objects $a, b$ of $\mathcal C$, the map of
simplicial sets
\[ \Map_{\mathcal C}(a,b) \rightarrow \Map_{\mathcal
D}(fa,fb) \] is a weak equivalence.

\item The induced functor on component categories $\pi_0f:\pi_0
\mathcal C \rightarrow \pi_0 \mathcal D$ is an equivalence of
categories.
\end{enumerate}
A Quillen equivalence between model categories then induces a
DK-equivalence between their simplicial localizations.

More generally, if one is not concerned with set-theoretic issues,
we can take the simplicial localization of any category with weak
equivalences.  Since, in a fairly natural sense, a ``homotopy
theory" is really some category with ``weak equivalences" that we
would like to invert, a homotopy theory gives rise to a simplicial
category.  In fact, the converse is also true: given any
simplicial category, it is, up to DK-equivalence, the simplicial
localization of some category with weak equivalences
\cite[2.5]{dkdiag}. Thus, the study of simplicial categories is
really the study of homotopy theories.

A first approach to applying the techniques of homotopy theory to
a category of simplicial categories itself was first given by
Dwyer and Kan \cite{dksl}.  In this paper, they define a model
category structure on the category of simplicial categories with a
fixed set $\mathcal O$ of objects.  The idea was then proposed
that the homotopy theory of (all) simplicial categories was
essentially the ``homotopy theory of homotopy theories."  Dwyer
and Spalinski mention this concept at the end of their survey
paper \cite{ds}, and the idea was further explored by Rezk
\cite{rezk}, whose ideas we will return to in the next section.
The author then showed in \cite{simpcat} that the category of all
small simplicial categories with the DK-equivalences has a model
category structure, thus formalizing the idea.

In order to define the fibrations in this model structure, we need
the following notion.  If $\mathcal C$ is a simplicial category
and $x$ and $y$ are objects of $\mathcal C$, a morphism $e \in
\Map_\mathcal C(x,y)_0$ is a \emph{homotopy equivalence} if the
image of $e$ in $\pi_0 \mathcal C$ is an isomorphism.

\begin{theorem} \cite[1.1]{simpcat}
There is a model category structure on the category $\mathcal{SC}$
of small simplicial categories defined by the following three
classes of morphisms:

\begin{enumerate}
\item The weak equivalences are the DK-equivalences.

\item The fibrations are the maps $f \colon \mathcal C \rightarrow
\mathcal D$ satisfying the following two conditions:
\begin{itemize}
\item For any objects $x$ and $y$ in $\mathcal C$, the map
\[ \Map_\mathcal C (x,y) \rightarrow \Map_\mathcal D (fx,fy) \]
is a fibration of simplicial sets.

\item For any object $x_1$ in $\mathcal C$, $y$ in $\mathcal D$,
and homotopy equivalence $e \colon fx_1 \rightarrow y$ in
$\mathcal D$, there is an object $x_2$ in $\mathcal C$ and
homotopy equivalence $d \colon x_1 \rightarrow x_2$ in $\mathcal
C$ such that $fd=e$.
\end{itemize}

\item The cofibrations are the maps which have the left lifting
property with respect to the maps which are fibrations and weak
equivalences.
\end{enumerate}
\end{theorem}

The advantage of this model category is that its objects are
fairly straightforward.  As mentioned above, there is a reasonable
argument for saying that simplicial categories really are homotopy
theories.  The disadvantage here lies in the weak equivalences, in
that they are difficult to identify.  Thus, it was natural to look
for a model with nicer weak equivalences.

\section{Complete Segal spaces}

Complete Segal spaces are probably the most complicated objects to
define of the four models described in this paper, but from the
point of view of homotopy theory, they might be the easiest to use
because the corresponding model structure gives what Dugger calls
a presentation for the homotopy theory \cite{dugger}. They are
defined by Rezk \cite{rezk} whose purpose was explicitly to find a
nice model for the homotopy theory of homotopy theories.

A complete Segal space is first a simplicial space.  It should be
noted that we require that certain of our objects be fibrant in
the Reedy model structure on the category of simplicial spaces
\cite{reedy}.  This structure is defined by levelwise weak
equivalences and cofibrations, but its importance here is that
several of our constructions will be homotopy invariant because
the objects involved satisfy this condition.

We begin by defining Segal spaces, for which we need the Segal map
assigned to a simplicial space. As one might guess from its name,
the Segal map is first defined by Segal in his work with
$\Gamma$-spaces \cite{segal}.  Let $\alpha^i:[1] \rightarrow [k]$
be the map in ${\bf \Delta}$ such that $\alpha^i(0)=i$ and
$\alpha^i(1)=i+1$, defined for each $0 \leq i \leq k-1$. We can
then define the dual maps $\alpha_i:[k]\rightarrow [1]$ in
$\Deltaop$.  For $k \geq 2$, the Segal map is defined to be the
map
\[ \varphi_k \colon X_k \rightarrow \underbrace{X_1 \times_{X_0} \cdots \times_{X_0} X_1}_k \]
induced by the maps
\[ X(\alpha_i) \colon X_k \rightarrow X_1. \]

\begin{definition} \cite[4.1]{rezk}
A Reedy fibrant simplicial space $W$ is a \emph{Segal space} if
for each $k \geq 2$ the map $\varphi_k$ is a weak equivalence of
simplicial sets. In other words, the Segal maps
\[ \varphi_k \colon W_k \rightarrow \underbrace{W_1 \times_{W_0} \cdots
\times_{W_0} W_1}_k \] are weak equivalences for all $k \geq 2$.
\end{definition}

A nice property of Segal spaces is the fact that they can be
regarded as analogous to simplicial categories, in that we can
define their ``objects" and ``morphisms" in a meaningful way.
Given a Segal space $W$, its set of objects, denoted $\ob(W)$, is
the set of 0-simplices of the space $W_0$, namely, the set
$W_{0,0}$. Given any two objects $x,y$ in $\ob (W)$, the mapping
space $\map_W(x,y)$ is the fiber of the map $(d_1, d_0) \colon W_1
\rightarrow W_0 \times W_0$ over $(x,y)$. Given a 0-simplex $x$ of
$W_0$, we denote by $\text{id}_x$ the image of the degeneracy map
$s_0 \colon W_0 \rightarrow W_1$. We say that two 0-simplices of
$\map_W(x,y)$, say $f$ and $g$, are homotopic, denoted $f \sim g$,
if they lie in the same component of the simplicial set
$\map_W(x,y)$.

Given $f \in \map_W (x,y)_0$ and $g \in \map_W(y,z)_0$, there is a
composite $g \circ f \in \map_W(x,z)_0$, and this notion of
composition is associative up to homotopy.  The homotopy category
$\Ho(W)$ of $W$, then, has as objects the set $\ob (W)$ and as
morphisms between any two objects $x$ and $y$, the set
$\map_{\Ho(W)}(x,y) = \pi_0 \map_W(x,y)$.

Finally, a map $g$ in $\map_W(x,y)_0$ is a homotopy equivalence if
there exist maps $f,h \in \map_W(y,x)_0$ such that $g \circ f \sim
\text{id}_y$ and $h \circ g \sim \text{id}_x$. Any map in the same
component as a homotopy equivalence is itself a homotopy
equivalence \cite[5.8]{rezk}.  Therefore we can define the space
$W_{\hoequiv}$ to be the subspace of $W_1$ given by the components
whose zero-simplices are homotopy equivalences.

We then note that the degeneracy map $s_0 \colon W_0 \rightarrow
W_1$ factors through $W_\hoequiv$ since for any object $x$ the map
$s_0(x)= \text{id}_x$ is a homotopy equivalence.  Therefore, we
have the following definition:

\begin{definition} \cite[\S 6]{rezk}
A \emph{complete Segal space} is a Segal space $W$ for which the
map $s_0 \colon W_0 \rightarrow W_\hoequiv$ is a weak equivalence
of simplicial sets.
\end{definition}

We can now consider some particular kinds of maps between Segal
spaces.  Note that, as the name suggests, these maps are very
similar in spirit to the weak equivalences in $\mathcal{SC}$.

\begin{definition}
A map $f \colon U \rightarrow V$ of Segal spaces is a
\emph{DK-equivalence} if
\begin{enumerate}
\item for any pair of objects $x,y \in U_0$, the induced map
\[ \map_U(x,y) \rightarrow \map_V(fx,fy)\] is a weak equivalence of
simplicial sets, and

\item the induced map $\Ho (f) \colon \Ho(U) \rightarrow \Ho(V)$
is an equivalence of categories.
\end{enumerate}
\end{definition}

We are now able to describe the important features of the complete
Segal space model category structure.

\begin{theorem}\cite[7.2, 7.7]{rezk}
There is a model structure $\css$ on the category of simplicial
spaces such that

\begin{enumerate}
\item The weak equivalences between Segal spaces are the
DK-equivalences.

\item The cofibrations are the monomorphisms.

\item The fibrant objects are the complete Segal spaces.
\end{enumerate}
\end{theorem}

What makes the model category $\css$ so nice to work with is the
fact that the weak equivalences between the fibrant objects, the
complete Segal spaces, are easy to identify.

\begin{prop} \cite[7.6]{rezk} \label{DKReedy}
A map $f \colon U \rightarrow V$ between complete Segal spaces is
a DK-equivalence if and only if it is a levelwise weak
equivalence.
\end{prop}

To avoid further technical detail, we have not defined what a
general weak equivalence is in $\css$, but the interested reader
can find it in Rezk's paper \cite[\S 7]{rezk}.  The important
point is that, when working with the complete Segal spaces, the
weak equivalences are especially convenient.

\section{Segal categories}

We now turn to our third model, that of Segal categories.  These
are natural generalizations of simplicial categories, in that they
can be regarded as simplicial categories with composition only
given up to homotopy. They first appear in the literature in a
paper of Dwyer, Kan, and Smith \cite{dks}, where they are called
special $\Deltaop$-diagrams of simplicial sets.  In particular,
Segal categories are again a kind of simplicial space.

We begin with the definition of a Segal precategory.

\begin{definition}
A \emph{Segal precategory} is a simplicial space $X$ such that
$X_0$ is a discrete simplicial set.
\end{definition}

As with the Segal spaces in the previous section, we can use the
Segal maps to define Segal categories.

\begin{definition}
A \emph{Segal category} $X$ is a Segal precategory $X: \Deltaop
\rightarrow \SSets$ such that for each $k \geq 2$ the Segal map
\[ \varphi_k \colon X_k \rightarrow \underbrace{X_1 \times_{X_0} \cdots \times_{X_0} X_1}_k \]
is a weak equivalence of simplicial sets.
\end{definition}

A model category structure $\Secat_c$ for Segal categories is
given by Hirschowitz and Simpson \cite{hs}.  In fact, they
generalize the definition to that of a Segal $n$-category and give
a model structure for Segal $n$-categories for any $n \geq 1$. The
idea behind this generalization is used for both the Simpson and
Tamsamani definitions of weak $n$-category \cite{simpson},
\cite{tam}.

The author gives a new proof of this model structure, just for the
case of Segal categories, from which it is easier to characterize
the fibrant objects \cite{thesis}.  It should be noted that, as in
the case of $\css$, this model structure is actually defined on
the larger category of Segal precategories.  However, the
fibrant-cofibrant objects are Segal categories \cite{fibrant}.

To define the weak equivalences in $\Secat_c$, we first note that
there is a functor $L_c$ assigning to every Segal precategory a
Segal category \cite[\S 5]{thesis}.  Then, if we are working with
a Segal category $X$, we can define its ``objects," "mapping
spaces," and ``homotopy category" just as we did for a Segal
space.

\begin{theorem}
There is a cofibrantly generated model category structure
$\Secat_c$ on the category of Segal precategories with the
following weak equivalences, fibrations, and cofibrations.
\begin{itemize}
\item Weak equivalences are the maps $f \colon X \rightarrow Y$
such that the induced map $\map_{L_cX}(x,y) \rightarrow
\map_{L_cY}(fx,fy)$ is a weak equivalence of simplicial sets for
any $x, y \in X_0$ and the map $\Ho(L_cX) \rightarrow \Ho(L_cY)$
is an equivalence of categories.

\item Cofibrations are the monomorphisms.  (In particular, every
Segal precategory is cofibrant.)

\item Fibrations are the maps with the right lifting property with
respect to the maps which are both cofibrations and weak
equivalences.
\end{itemize}
\end{theorem}

It should not be too surprising that the weak equivalences in this
case are again called DK-equivalences, since the same idea
underlies the definition in each of the three categories we have
considered.

Furthermore, there is also second model structure $\Secat_f$ on
this same category, with the same weak equivalences but different
fibrations and cofibrations.  Thus each leads to the same homotopy
theory, and in fact they are Quillen equivalent, but the slight
difference between the two is key in comparing Segal categories
with the other models. We will not define the fibrations and
cofibrations in $\Secat_f$ here, as they are technical and
unenlightening in themselves, but we refer the interested reader
to \cite[\S 7]{thesis} for the details.  As the subscript
suggests, the initial motivation was to find a model structure
with the same weak equivalences but in which the fibrations,
rather than the cofibrations, were given by levelwise fibrations
of simplicial sets.  As it turns out, such a description does not
work, but one is not far off thinking of the fibrations as being
levelwise.

\section{Quasi-categories}

Quasi-categories are perhaps the most mysterious, as far as why
they are equivalent to the others. While each of the other models
consists of simplicial spaces, or objects easily related to
simplicial spaces, quasi-categories are simplicial sets, and thus
simpler than any of the others. Thus, they may also provide a good
model to use when we actually want to compute something for a
given homotopy theory.  They were first defined by Boardman and
Vogt \cite{bv} and are sometimes called weak Kan complexes.

Recall that in the category of simplicial sets we have several
particularly important objects.  For each $n \geq 0$, there is the
$n$-simplex $\Delta[n]$ and its boundary $\dot \Delta [n]$.  If we
remove the $k$th face of $\dot \Delta [n]$, we get the simplicial
set denoted $V[n,k]$.  Given any simplicial set $X$, a horn in $X$
is a map $V[n,k] \rightarrow X$.  A Kan complex is a simplicial
set such that every horn factors through the inclusion map $V[n,k]
\rightarrow \Delta [n]$.  A quasi-category is then a simplicial
set $X$ such that every horn $V[n,k] \rightarrow X$ factors
through $\Delta [n]$ for each $0<k<n$.  Such horns are called the
inner horns.  Further details on quasi-categories can be found in
Joyal's papers \cite{joyal1} and \cite{joyal2}.

Like the cases for Segal categories and complete Segal spaces, the
model structure $\Qcat$ for quasi-categories is defined on a
larger category, in this case the category of simplicial sets.  To
define the weak equivalences, we need some definitions.

First, consider the nerve functor $N: \Cat \rightarrow \SSets$,
where $\Cat$ denotes the category of small categories.  This
functor has a left adjoint $\tau_1:\SSets \rightarrow \Cat$. Given
a simplicial set $X$, the category $\tau_1(X)$ is called its
fundamental category.  We can then define a functor $\tau_0:\SSets
\rightarrow \Sets$, where $\tau_0(X)$ is the set of isomorphism
classes of objects of the category $\tau_1(X)$.  Now, if $X^Y$
denotes the simplicial set of maps $Y \rightarrow X$, for any pair
$(X,Y)$ of simplicial sets we define $\tau_0(Y,X)= \tau_0(X^Y)$. A
weak categorical equivalence is a map $A \rightarrow B$ of
simplicial sets such that the induced map $\tau_0(B,X) \rightarrow
\tau_0(A,X)$ is an isomorphism of sets for any quasi-category $X$.

\begin{theorem} \cite{joyal2}
There is a model category structure $\Qcat$ on the category of
simplicial sets in which the weak equivalences are the weak
categorical equivalences and the cofibrations are the
monomorphisms. The fibrant objects of $\Qcat$ are the
quasi-categories.
\end{theorem}

\section{Quillen equivalences}

The origins of the comparisons between these various structures
seem to be in various places.  The question of simplicial
categories and Segal categories is a fairly natural one, since a
Segal category is essentially a simplicial category up to
homotopy.  It is addressed partially by Dwyer, Kan, and Smith
\cite{dks}, but they do not give a Quillen equivalence, partly
because their work predates both model structures by several
years.  Schw\"anzl and Vogt also address this question, using
topological rather than simplicial categories \cite{sv}.

Rezk defines complete Segal spaces with the comparison with
simplicial categories in mind \cite{rezk}.  While his functor from
simplicial categories to complete Segal spaces naturally factors
through Segal categories, he does not mention this fact as such.
Initially, there did not seem to be a need to bring in the Segal
categories from this point of view, but further investigation led
to skepticism that his functor had the necessary adjoint to give a
Quillen equivalence.

To\"{e}n mentions all four models and conjectures the
relationships between them in \cite{toen1}. As mentioned in the
introduction, some but not all of these functors are the ones used
in the known proofs.  To\"{e}n further axiomatizes the notion of
``theory of $(\infty,1)$-categories" in \cite{toen2}.  He gives
six axioms for a model structure to satisfy in order to be such a
theory, and he shows that any such model category structure is
Quillen equivalent to Rezk's complete Segal space model structure.
As far as the author knows, these axioms have not been verified
for the other three models given in this paper, but it seems
likely that they should hold.

Some of the comparisons are also mentioned in work by Simpson, who
sketches an argument for comparing the Segal categories and
complete Segal spaces \cite{simpson2}.

The author showed in \cite{thesis} that simplicial categories are
equivalent to Segal categories, which are in turn equivalent to
the complete Segal spaces.  However, the adjoint pairs go in
opposite directions and therefore cannot be composed into a single
Quillen equivalence.  It is still unknown, as far as we know,
whether there is a direct Quillen equivalence.  Further work by
Joyal and Tierney has shown that there are Quillen equivalences
between quasi-categories and each of the other three models.

We now look in more detail at these Quillen equivalences.  Let us
begin with the Segal categories and complete Segal spaces.  Since
the underlying category of $\css$ is the category of simplicial
spaces and the underlying category of $\Secat_c$ is the category
of simplicial spaces with a discrete space in degree zero, there
is an inclusion functor $I \colon \Secat_c \rightarrow \css$. This
functor has a right adjoint $R \colon \css \rightarrow \Secat_c$
which acts as a discretization functor.  In particular, if it is
applied to a complete Segal space $W$, the result is a Segal
category which is DK-equivalent to it.

\begin{theorem} \cite[6.3]{thesis}
The adjoint pair
\[ \xymatrix@1{I \colon \Secat_c \ar@<.5ex>[r] & \css :R \ar@<.5ex>[l]} \] is
a Quillen equivalence.
\end{theorem}

Then, as we mentioned in the section on Segal categories, the two
model structures $\Secat_c$ and $\Secat_f$ are Quillen equivalent.

\begin{theorem} \cite[7.5]{thesis}
The identity functor induces a Quillen equivalence
\[ \xymatrix@1{id \colon \Secat_f \ar@<.5ex>[r] & \Secat_c:id. \ar@<.5ex>[l]} \]
\end{theorem}

Turning to simplicial categories, the nerve functor $N:
\mathcal{SC} \rightarrow \Secat_f$ has a left adjoint $F$, which
can be considered a rigidification functor.

\begin{theorem} \cite[8.6]{thesis}
The adjoint pair
\[ \xymatrix@1{F \colon \Secat_f \ar@<.5ex>[r] & \mathcal{SC} :N \ar@<.5ex>[l]} \] is
a Quillen equivalence.
\end{theorem}

Turning to the quasi-categories, Joyal and Tierney have shown that
there are in fact two different Quillen equivalences between
$\Qcat$ and $\css$.  For the first of these equivalences, the map
$i^*_1$, which associates to a complete Segal space $W$ the
simplicial set $W_{*0}$, has a left adjoint $p^*_1$.

\begin{theorem} \cite{jt} \label{qcat1}
The adjoint pair of functors
\[ \xymatrix@1{p^*_1 \colon \Qcat \ar@<.5ex>[r] & \css :i^*_1 \ar@<.5ex>[l]} \] is
a Quillen equivalence.
\end{theorem}

The second Quillen equivalence between these two model categories
is given by a total simplicial set functor $t_!:\css \rightarrow
\Qcat$ and its right adjoint $t^!$.

\begin{theorem} \cite{jt}
The adjoint pair
\[ \xymatrix@1{t_! \colon  \css \ar@<.5ex>[r] & \Qcat :t^! \ar@<.5ex>[l]}
\] is a Quillen equivalence.
\end{theorem}

Even one of these Quillen equivalences would be sufficient to show
that all four of our model categories are equivalent to one
another, but, interestingly, Joyal and Tierney go on to prove that
there are also two different Quillen equivalences directly between
$\Qcat$ and $\Secat_c$.  The first of these functors is analogous
to the pair given in Theorem \ref{qcat1}; the right adjoint
functor $j^* \colon \Secat_c \rightarrow \Qcat$ assigns to a Segal
precategory $X$ the simplicial set $X_{*0}$.  Its left adjoint is
denoted $q^*$.

\begin{theorem} \cite{jt}
The adjoint pair
\[ \xymatrix@1{q^* \colon \Qcat \ar@<.5ex>[r] & \Secat_c :j^* \ar@<.5ex>[l]} \]
is a Quillen equivalence.
\end{theorem}

The second Quillen equivalence between these two model categories
is given by the map $d^* \colon \Secat_c \rightarrow \Qcat$, which
assigns to a Segal precategory its diagonal, and its right adjoint
$d_*$.

\begin{theorem} \cite{jt}
The adjoint pair
\[ \xymatrix@1{d^* \colon \Secat_c \ar@<.5ex>[r] & \Qcat :d_* \ar@<.5ex>[l]}
\] is a Quillen equivalence.
\end{theorem}

Finally, Joyal has also related the quasi-categories to the
simplicial categories directly.  There is the coherent nerve
functor $\widetilde{N} \colon \mathcal{SC} \rightarrow \Qcat$,
first defined by Cordier and Porter \cite{cp}. Given a simplicial
category $X$ and the simplicial resolution $C_*[n]$ of the
category $[n]=(0 \rightarrow \cdots \rightarrow n)$, the coherent
nerve $\widetilde{N}(X)$ is defined by
\[ \widetilde{N}(X)_n = \Hom_{\mathcal SC}(C_*[n],X). \]  This
functor has a left adjoint $J \colon \Qcat \rightarrow
\mathcal{SC}$.

\begin{theorem} \cite{joyal3}
The adjoint pair
\[ \xymatrix@1{J \colon \Qcat \ar@<.5ex>[r] & \mathcal{SC} :\widetilde{N} \ar@<.5ex>[l]}
\] is a Quillen equivalence.
\end{theorem}

Thus, we have the following diagram of Quillen equivalences of
model categories:
\[ \xymatrix{\mathcal{SC} \ar@<-.5ex>[r] \ar@<-.5ex>[drr] & \Secat_f
\ar@<.5ex>[r] \ar@<-.5ex>[l]  & \Secat_c \ar@<.5ex>[r]
\ar@<.5ex>[l] \ar@{<->}[d] & \css \ar@<.5ex>[l] \ar@{<->}[ld] \\
&& \Qcat \ar@{<->}[u] \ar@{<->}[ur] \ar@<-.5ex>[ull]} \]

The single double-headed arrows indicate that in these cases
either direction can be chosen to be a left (or right) adjoint,
depending on which Quillen equivalence is used.

\end{document}